\documentclass[a4paper]{article}

\author{Joaquim Ro\'e \\
\small{Departament de Matem\`{a}tiques, Universitat Auton\`{o}ma de
Barcelona,}\\ \small{Edifici C, 08193-Bellaterra
  (Barcelona), Spain. e-mail: jroe@mat.uab.es}}

\usepackage[centertags]{amsmath}
\usepackage{amsthm}
\usepackage{amssymb}

\theoremstyle{definition}
\newtheorem{Def}{Definition}%[section]
\theoremstyle{plain}

\newtheorem{Cor}[Def]{Corollary}
\newtheorem{Pro}[Def]{Proposition}
\newtheorem{Teo}[Def]{Theorem}
\newtheorem{Con}[Def]{Conjecture}
\theoremstyle{remark}
\newtheorem{Rem}[Def]{Remark}

\newcommand{\mult}{\operatorname{mult}}

\newcommand{\A}{{\mathbb A}}

\newcommand{\M}{\mathfrak m}
\newcommand{\R}{\mathbb{R}}
\newcommand{\Z}{\mathbb{Z}}
\newcommand{\C}{\mathbb{C}}
\newcommand{\m}{{\bf m}}
\newcommand{\x}{{\bf x}}

\renewcommand{\P}{\mathbb{P}}
\renewcommand{\H}{{\cal H}}
\renewcommand{\O}{{\cal O}}

\begin{document}

\title{On Nagata's conjecture}
\maketitle

\begin{abstract}
T. Szemberg proposed in 2001 a generalization to arbitrary varieties of M.
Nagata's 1959 open conjecture, which claims that the Seshadri constant of $r\ge
9$ very general points of the projective plane is maximal. Here we prove that
Nagata's original cojecture implies Szemberg's  for all smooth surfaces $X$ with
an ample divisor $L$ generating $NS(X)$ and such that $L^2$ is a square.

More generally, we prove the inequality
$$\varepsilon_{n-1}(L,r)\ge\varepsilon_{n-1}(L,1)\varepsilon_{n-1}(\O_{\P^n}(1),r),$$
where $\varepsilon_{n-1}(L,r)$ stands for the $(n-1)$-dimensional Seshadri
constant of the ample divisor $L$ at $r$ very general points of a normal
projective variety $X$, and $n=\dim X$.
\end{abstract}

\section{Introduction}

Let $X$ be a normal projective variety of dimension $n$ over an algebraically
closed field $k$, and $L$ an ample divisor. Given $r$ points $p_1, \dots, p_r \in
X$ and an integer $1 \le d \le n$, the $d$-dimensional Seshadri constant of $L$
at the points $p_1, \dots, p_r$ is the real number
$$
\varepsilon_d(L,p_1, \dots, p_r)= \sqrt[d]{ \inf_Z \left\{\frac{L^d \cdot Z}{\sum
\mult_{p_i}Z}\right\}},$$ where $Z$ runs over all positive $d$-dimensional cycles
(see \cite[1.1]{Sze01}, or \cite{Dem92} for the original definition). As $L$ is
ample, we have $L^d \cdot Z>0$ for all $Z$, so $\frac{L^d \cdot Z}{\sum
\mult_{p_i}Z} \in \R_+ \cup \infty$. Moreover, there exist $Z$ which contain some
point $p_i$, and therefore the Seshadri constant is indeed a finite real number.

Most work on Seshadri constants deals with the $d=1$ case, and usually one writes
$\varepsilon(L,p_1, \dots, p_r)= \varepsilon_1(L,p_1, \dots, p_r)$; we shall be
concerned here with the codimension 1 case $(d=n-1)$. Also, we use the shorthand
notation $\varepsilon_d(L,r)=\varepsilon_d(L,p_1, \dots, p_r)$ for very general
points $p_1, \dots, p_r$ (i.e., in the intersection of countably many Zariski
open subsets of $X^{r}$) which is the case we are mostly interested in.
%As Demailly
%pointed out in \cite[Remark 6.7]{Dem92}, one has $\varepsilon_d(L,p_1, \dots,
%p_r) \ge \varepsilon(L,p_1, \dots, p_r)$ for every $d$.

In connection with his solution to the fourteenth problem of Hilbert, Nagata
posed in \cite{Nag59} (in different terminology) the following conjecture
concerning Seshadri constants of the plane:

\begin{Con}[Nagata]
If $r\ge 9$, then
$$\varepsilon(\O_{P^2}(1),r)=1/\sqrt{r}.$$
\end{Con}

If $r=s^2$ is a square, then it is not hard to prove that the conjecture is true,
and in fact Nagata proved a slightly stronger result in this case. In a variety
of dimension $n$ it is also not difficult to prove that
$$\varepsilon(L,p_1, \dots, p_r)\le \sqrt[n]{\frac{L^n}{r}}$$
for every set of $r$ points (see \cite[Remark 1]{Ste98}, for example), so
Nagata's conjecture claims that the Seshadri constant of a very general set of
$r\ge 9$ points in the plane is maximal. All available information on Seshadri
constants (see \cite{Xu95}, \cite{Laz93}, \cite{Kuc96a}, \cite{Bir98},
\cite{Bau99}, \cite{ST?2}, \cite{Har?2} for the case of surfaces, \cite{AH00},
\cite{Chu78}, \cite{HK??} for dimension $n>2$) suggests that, in fact, in an
arbitrary variety, for $r$ large enough, the Seshadri constant of $r$ very
general points is maximal. This led Szemberg to propose in \cite{Sze01} the
following generalization:

\begin{Con}[Nagata-Szemberg]
Given a smooth variety of dimension $n$ and $L$ an ample divisor on $X$, there
exists a number $r_0=r_0(X,L)$ such that for every $r\ge r_0$
$$\varepsilon(L,r)=\sqrt[n]{\frac{L^n}{r}}.$$
\end{Con}

This note is devoted to the following result, that gives a lower bound for
$(n-1)$-dimensional Seshadri constants of $r$ very general points in a variety,
relating them to the analogous constants of one point in the same variety and of
$r$ points in projective space:

\begin{Teo}
\label{main} Let $X$ be a normal projective variety of dimension
$n \ge 2$, $L$ an ample divisor. Then for every smooth point $p
\in X$ and $r\ge 1$,
$$\varepsilon_{n-1}(L,r)\ge\varepsilon_{n-1}(L,p)\varepsilon_{n-1}(\O_{\P^n}(1),r).$$
\end{Teo}

Combining it with known results on the value of the Seshadri constants in
projective space, Theorem \ref{main} implies more explicit relations between
$r$-point and 1-point Seshadri constants. The consequences support the
Nagata-Szemberg conjecture, especially in the case of surfaces (note that for
surfaces, $(n-1)$-dimensional Seshadri constants are the usual Seshadri
constants).

\begin{Cor}
Suppose $r=s^n$ for some integer $s$. Then for every smooth point $p \in X$ and
$r\ge 1$,
$$\varepsilon_{n-1}(L,r)\ge\varepsilon_{n-1}(L,p)/s.$$
\end{Cor}
\begin{proof}
Use G. V. Choodnovsky's result \cite{Chu78} that
$\varepsilon_{n-1}(\O_{\P^n}(1),s^n)=1/s$.
\end{proof}

If $r$ is not the $n$-th power of an integer, then we do not know the exact value
of $\varepsilon_{n-1}(\O_{\P^n}(1),r)$ (there is a conjecture similar to
Nagata's, posed by Choodnovsky in the same paper \cite{Chu78}, and by A.
Iarrobino in \cite{Iar97}). However, B. Harbourne pointed out that,
 using results of J. Alexander and A.
Hirschowitz and of M. Hochster and C. Huneke, one can prove an asymptotically
optimal bound for $\varepsilon_{n-1}(\O_{\P^n}(1),r)$. Combining it
with Theorem \ref{main}, in section \ref{proofs} we prove the following
asymptotic bound for $\varepsilon_{n-1}(L,r)$ that depends only on the
$(n-1)$-dimensional Seshadri constant of $L$ at a smooth point $p\in
X$:

\begin{Cor}
\label{ahhh} For every $\epsilon>0$ there exists an integer $s$, depending only
on $\epsilon$ and $n$, such that for every smooth point $p\in X$ and $r\ge s$,
$$\varepsilon_{n-1}(L,r){\sqrt[n]{r}}\ge \varepsilon_{n-1}(L,p)- \epsilon.$$
\end{Cor}

\begin{Cor}
If $X$ is a normal projective surface, $L$ an ample divisor on $X$, then for
every smooth point $p \in X$ and $r\ge 9$, Nagata's conjecture implies that
$$\varepsilon(L,r)\ge\frac{\varepsilon(L,p)}{\sqrt{r}}.$$
\end{Cor}

Observe that this tells us that Nagata's conjecture (on the plane) implies the
Nagata-Szemberg conjecture on a large family of surfaces. Indeed, the obtained
bound is equal to $\varepsilon(L,p)/\sqrt{L^2}$ times the conjectured value of
$\varepsilon(L,r)$, so we have the following:
\begin{Cor}
If $X$ is a normal projective surface, $L$ an ample divisor on $X$ and $p\in X$
is a smooth point such that $\varepsilon(L,p)=\sqrt{L^2}$, then Nagata's
conjecture implies the Nagata-Szemberg conjecture on $(X,L)$ with $r_0(X,L)\le
9$.
\end{Cor}

In particular, we can apply this to complex surfaces with Picard number 1, using
A. Steffens' result \cite[Proposition 1]{Ste98}, which says that if $L$ is an
ample generator of $NS(X)$ then $\varepsilon(L,p)\ge \lfloor \sqrt{L^2}\rfloor$
for very general points:
\begin{Cor}
Let $X$ be a smooth projective surface defined over $\C$, $L$ an ample generator
of $NS(X)$ and assume $L^2=d^2$ is a square. Then Nagata's conjecture implies the
Nagata-Szemberg conjecture on $X$ with $r_0(X,L)\le 9$.
\end{Cor}

In particular, the Seshadri constant of $r\ge 9$ very general points on a complex
surface $X$ with Picard number equal to 1, is maximal if both $L^2$ and the
number of points $r$ are squares. This can be compared to B. Harbourne's result
\cite[I.2.b(i)]{Har?2} (over an arbitrary base field and with no assumption on
the Picard number) that the Seshadri constant is maximal whenever $L$ is very
ample, $rL^2$ is a square, and $r\ge L^2$.

Also, known bounds approximating Nagata's conjecture give new bounds on surfaces;
for instance, H. Tutaj-Gasi\'nska's bound in \cite{Tut??} showing that
$\varepsilon(\O_{\P^2_{\C}}(1),r)\ge\frac{1}{\sqrt{r+\frac{1}{12}}}$ gives the
following:

\begin{Cor}
\label{tg} If $X$ is a normal projective surface defined over $\C$, then for
every smooth point $p \in X$ and $r>9$,
$$\varepsilon(L,r)\ge\frac{\varepsilon(L,p)}{\sqrt{r+\frac{1}{12}}}.$$
\end{Cor}

In a similar vein, Harbourne's bounds on Seshadri constants of $\P^2$ in
\cite{Har01} and \cite{Har?2} imply that
\begin{Cor}
\label{ha} Let $X$ be a normal projective surface, $p \in X$ a smooth point and
$r\ge 1$. Then for every pair of integers $1 \le s \le r$, $1 \le d$, it holds
$$\varepsilon(L,r)\ge
\begin{cases}
\frac{s}{rd}\varepsilon(L,p) & \text{ if }s^2 \le rd^2, \\
\frac{d}{s}\varepsilon(L,p) & \text{ if }s^2 \ge rd^2.
\end{cases}
$$
\end{Cor}
It should be mentioned, however, that both Corollaries \ref{tg} and \ref{ha} are
usually weaker (but stronger for some surfaces and numbers of points) than
Harbourne's results on algebraic surfaces of \cite{Har?2}, where he gets the
bounds
$$\varepsilon(L,r)\ge
\begin{cases}
\frac{s}{rd} & \text{ if }s^2 \le rd^2L^2, \\
\frac{dL^2}{s} & \text{ if }s^2 \ge rd^2L^2
\end{cases}
$$
for very ample $L$, assuming moreover that $r \ge L^2$.

The proof of Theorem \ref{main} is based on the idea, due to L. \'Evain (see
\cite{Eva97}) that $r$ point Seshadri constants of the plane can be computed by
means of homotetic collisions of fat points. For convenience of the exposition we
express this, generalized to $n$-dimensional projective space, in terms of the
order of a nonreduced curve singularity, rather than collisions. Then, we observe
that it is enough to know the formal germ of such a singularity, and the fact
that the completion of the local ring at a smooth point of a variety is a ring of
formal power series that only depends on the dimension of the variety, to obtain
the bound.

It might be interesting to note that the Viro method developed by E. Shustin in
\cite{Shu98} can also be used to relate the existence of singular curves on
$\P^2$ with the existence of singular curves in algebraic surfaces, see for
instance \cite[\S 5]{KT02} or \cite[3.A]{Bir98}.

We thank B. Harbourne for many valuable comments, which largely
improved the paper.

%\section{Notation and setup}

\section{Singularities of arrangements of multiple lines}
\label{proofs}

To lighten notations for rings of polynomials and of power series, we write
$\x=(x_0,$ \dots, $x_n)$ for a collection of variables, so $k[\x]$ and $k[[\x]]$
denote $k[x_0,\dots,x_n]$ and $k[[x_0,\dots,x_n]]$, respectively. Also, $\M$ and
$\hat \M$ will be the maximal ideals generated by $x_0,\dots,x_n$ in $k[\x]$ and
$k[[\x]]$ respectively, and for every point $p=[\xi_0:\dots:\xi_n]$ in projective
$n$-space, we use the notation $I_p$ (respectively $\hat I_p$) for the ideal
generated by the $2\times 2$ minors of the matrix
%$(\{\xi_jx_i-\xi_ix_j\}_{i\ne j})$
$$
\begin{pmatrix}
\xi_0 & \dots & \xi_n \\
x_0 & \dots & x_n
\end{pmatrix}
$$
in $k[\x]$ (respectively in $k[[\x]]$).

Given distinct points $p_1, \dots, p_r \in \P^n$ and $\m=(m_1, \dots, m_r) \in
\Z^r_{\ge 0}$, we define $\alpha_{\m}(p_1, \dots, p_r)$ to be the minimal degree
of a homogeneous polynomial vanishing to order $m_i$ at $p_i$. As
$I_p$ is homogeneous for all $p$, this
number coincides with the maximal integer $\alpha$ such that
$$
I=\bigcap_{i=1}^r I_{p_i}^{m_i}\subset \M^\alpha,
$$
or equivalently, such that $\hat I= \bigcap \hat
I_{p_i}^{m_i}\subset \hat \M^\alpha$. In other words, $\alpha_{\m}(p_1,
\dots, p_r)$ is the order at the origin of $\A^{n+1}$ of the arrangement of
multiple lines defined by $I$ (which is the affine cone over the fat point scheme
consisting of the points $p_i$ with multiplicities $m_i$).

\begin{Rem} \label{alpha} The definition of $\alpha$ and the $(n-1)$-dimensional Seshadri
constants immediately give that $\forall \m$,
$$\alpha_\m(p_1, \dots, p_r)\ge \left(\varepsilon_{n-1}(\O_{\P^n}(1),p_1, \dots,
p_r)\right)^{n-1}\sum_{i=1}^r m_i .$$
\end{Rem}

Let $X$ be a variety of dimension $n\ge 2$, $q\in X$ a smooth point, and fix
uniformizing parameters $x_1, \dots, x_n$ in some neighborhood $V$ centered at
$q$. So (see \cite[\S III.6]{Mum99}, for example) the germs of $x_1, \dots,
x_n$ at $q$ generate the maximal ideal in the local ring $\O_{X,q}$ and the
morphism of $k$-algebras $k[\x]\hookrightarrow \O_X(V)$ determines an \'etale
morphism $\varphi:V \longrightarrow \A^n \cong T_q X$ (where $V\subset X$ is
open) and an isomorphism $k[[\x]] \longrightarrow \smash{\hat\O_{X,q}}$ .

To every $p=[\xi_0:\dots:\xi_n]\in \P^n$ such that $\xi_0 \ne 0$ and
$p\ne[1:0:\dots:0]$ (so that $p'=[\xi_1:\dots:\xi_n]\in \P^{n-1}$), we shall
assign an irreducible curve $C_p$, smooth at $q$, and a regular parameter $\bar
x_p\in O_{C_p,q}$. Define $C_p$ as the closure of the component through $q$ of
$\varphi^{-1}(L_{p'})$, where $L_{p'}\subset \A^{n}$ is the affine cone over
$p'$, i.e., the line through the origin in the direction determined by $p'$. As
$\varphi$ is \'etale, $C_p$ is smooth at $q$, and the ideal $I_{C_p} \subset
\O_{X,q}$ of its germ at $q$ is generated by the $2\times 2$ minors of the matrix
%$(\{\xi_jx_i-\xi_ix_j\}_{i\ne j})$
$$
\begin{pmatrix}
\xi_1 & \dots & \xi_n \\
x_1 & \dots & x_n
\end{pmatrix}.
$$
% The assumption
% $p \ne [1:0:\dots:0]$ guarantees that $C_p$ is indeed a curve, and $\xi_0 \ne 0$
% allows us to
Now consider $x_p=x_i \xi_i/\xi_0$ for some $\xi_i \ne 0$, $i \ge 1$.
It is easy to see that the restriction $\bar x_p$ of $x_p$ to $C_p$
does not depend on the choice of $i$, and that it is a uniformizing
parameter for the curve.
% It is easy to see that $\bar x_p \in \O_{C_p,q}=\O_{X,q}/I_{C_p}$
% does not depend on the choice of $i$, and that it is
% a regular parameter of $\O_{C_p,q}$.
Also, abusing notation,
in $\O_{C_p \times X,(q,q)}=\O_{C_p,q} \otimes \O_{X,q}$ we write
$x_p=\bar x_p \otimes 1$, $x_1 = 1 \otimes x_1$, \dots, and $x_n = 1
\otimes x_n$; then the ideal of the germ of the diagonal
$\Delta(C_p)\subset C_p \times X$ is generated by the $2\times 2$
minors of the matrix
$$
\begin{pmatrix}
\xi_0 & \xi_1 & \dots & \xi_n \\
x_p & x_1 & \dots & x_n
\end{pmatrix},
$$
in other words, $\Delta(C_p)$ is the closure of the
component through the origin
of $\psi^{-1}(L_p)$, where $L_p$ is the affine cone over $p \in \P^n$,
i.e., a line $L_p \subset \A^{n+1}=T_{(q,q)}(C_p \times X)$,
and $\psi$ is the \'etale morphism given by the parameters $x_p$, $x_1$,
\dots, $x_n$.

With these notations, Theorem \ref{main} follows from the more precise
proposition:

\begin{Pro}
\label{variety} Let $X$ be a variety of dimension $n\ge 2$, $L$ an ample divisor,
$q \in X$ a smooth point, $p_1, \dots, p_r \in \P^n  \setminus [1:0:\dots:0]$
distinct points not on the hyperplane $\xi_0=0$. Then for very general points
$q_k \in C_{p_k}$ it holds
$$\varepsilon_{n-1}(L, q_1, \dots, q_r)\ge \varepsilon_{n-1}(L,q)
\varepsilon_{n-1}(\O_{\P^n}(1),p_1, \dots, p_r).$$
\end{Pro}
\begin{proof}
First note that due to the semicontinuity of multiplicity (see \cite[\S 8]{LT73} or
\cite[\S 3]{Lip82}), for each component $\H$ of the Hilbert scheme of hypersurfaces
in $X$, and each system of multiplicities $\m$, the sets of points $(q_1, \dots,
q_r)$ such that there is $Y\in \H$ with multiplicity $\ge m_i$ at $q_i$ form a
Zariski-closed subset of $X^r$. Thus, it will be enough to prove that, given $\H$
and $\m$, the existence of $Y\in \H$ with multiplicity $\ge m_k$ at $q_k$ for
general $q_k \in C_{p_k}$ implies
$$Y\cdot L^{n-1}\ge \left(\varepsilon_{n-1}(L,q)
\varepsilon_{n-1}(\O_{\P^n}(1),p_1, \dots, p_r)\right)^{n-1}\sum_{k=1}^r m_k.$$

So, fix $\H$ and $\m$, and assume that for general points $q_k \in C_{p_k}$ there
is a hypersurface $Y \in \H$ going through $q_k$ with multiplicity at least
$m_k$.

In the local ring of $\prod C_{p_k}$ at $\Delta(q)=(q, \dots, q)$, the $1
\otimes \dots \otimes \bar x_{p_k} \otimes \dots \otimes 1$,
$k=1, \dots,r$, form a regular
system of parameters; abusing notation we call them simply $x_{p_k}$.
Let $\Gamma \subset \prod C_{p_k}$ be the irreducible curve defined locally by the
equations $x_{p_1}=\dots=x_{p_r}$, which is obviously smooth at
$\Delta(q)$, and admits $x_0:=\bar x_{p_1} \in \O_{\Gamma,\Delta(q)}$ as a
local parameter. For every $p_k=[\xi_{k,0}: \dots : \xi_{k,n}]$,
the ideal $J_{p_k} \subset \O_{\Gamma\times X,\Delta(q)}$ generated by
the $2\times 2$ minors of the matrix
$$
\begin{pmatrix}
\xi_{k,0} &  \dots & \xi_{k,n} \\
x_0 & \dots & x_n
\end{pmatrix},
$$
defines the germ of a curve $C_{p_k}'$ (the preimage of a line in $\A^{n+1}\cong
T_{\Delta(q)}(\Gamma \times X)$) whose projection to $X$ is exactly $C_{p_k}$;
more precisely, the fiber of $C_{p_k}'$ over $\gamma=(q_1, \dots, q_n) \in
\Gamma$ is $q_k \in C_{p_k}$.

Consider now the diagonals $\Delta_{i,j}=\{(q_1, \dots, q_r) \in \prod C_{p_k}|
q_i=q_j\}$, $\Delta=\bigcup \Delta_{i,j}$, and $U=\Gamma\setminus \Delta$. As the
points $p_1$, \dots, $p_r$ are distinct, $U$ is a nonempty open subset of
$\Gamma$. The assumption on the existence of $Y$ tells us that there is an
effective relative Weil divisor ${\mathbf Y}_U\subset U \times X$ (flat over $U$)
whose fiber over $(q_1, \dots, q_r)\in U$ belongs to $\H$ and has multiplicity at
least $m_i$ at $q_i$. By the smoothness of $\Gamma$ at $\Delta(q)$, ${\mathbf
Y}_U$ can be extended to a flat family ${\mathbf Y}\subset {U \cup \{\Delta(q)\}}
\times X$, and then the condition on the multiplicity of the fibers of $\mathbf
Y$ means that $\mathbf Y$ contains the arrangement of multiple curves whose germ
at $\Delta(q)$ is defined by the ideal
$$
J=\bigcap_{k=1}^r J_{p_k}^{m_k},
$$
and this implies that the fiber $Y_q$ of ${\mathbf Y}$ over $\Delta(q)$ has
multiplicity at least equal to the order of this arrangement at
$\Delta(q)$, i.e., at
least equal to the maximal integer $\alpha$ such that
$J\subset(x_0,\dots,x_n)^\alpha$. This can be computed equivalently as the order
of the completion $\hat J$ in $\hat \O_{\A^1\times X,(0,q)}\cong k[[\x]]$ but, by
construction of the curves $C_{p_k}'$, one has $\hat J=\hat I$ as defined above,
so this order is exactly $\alpha_\m(p_1, \dots, p_r)$.

We remark also that, by the smoothness of $X$ in a neighborhood of $q$, ${\mathbf
Y}$ is defined by a principal ideal at $q$, so we get a Weil divisor
$Y_q'\subseteq Y_q$ with
\begin{align*}
L^{n-1} \cdot Y_q' & = L^{n-1} \cdot Y, \\
\mult_q Y_q' & \ge \alpha_\m(p_1, \dots, p_r),
\end{align*}
which together with the definition of the 1-point Seshadri constant, gives
$L^{n-1} \cdot Y \ge \varepsilon_{n-1}(X,q)\alpha_\m(p_1, \dots, p_r)$, and then
it is enough to apply the bound of Remark \ref{alpha}.
\end{proof}

\begin{proof}[Proof of Corollary \ref{ahhh}]
By Theorem \ref{main} it is enough to see that
$$ \lim_{r\rightarrow \infty} \varepsilon_{n-1}(\O_{\P^n}(1),r)\sqrt[n]{r} \ge 1 $$
(which in fact means that one has an equality, the converse
inequality being well-known). More precisely, we shall prove that
given $k>0$ there exists $s_k=s_k(n)$ such that if $r \ge s_k$
then for all $\m=(m_1, \dots, m_r)$ and general points $p_1$,
\dots, $p_r$,
$$\alpha_{\m}(p_1, \dots,p_r)\ge \frac{\sum m_i
}{\sqrt[n]{r^{n-1}}}\cdot{\frac{k+1}{k+n}}.$$ So let $F$ be a
homogeneous polynomial defining a hypersurface of degree $d$ in
$\P^n$ which has multiplicity $m_i$ at $p_i$ for general points
$p_1$, \dots, $p_r$. Then, by the genericity of the points, for
every permutation $\sigma \in {\cal S}_r$ there is a polynomial
$F_\sigma$ which has multiplicity $m_i$ at the point
$p_{\sigma(i)}$. Therefore $G=\prod_{\sigma\in {\cal S}_r}
F_{\sigma}$ is a polynomial of degree $D=r!\,d$ which has (the
same) multiplicity $M=(r-1)!\sum m_i$ at $p_1$, \dots, $p_r$, and
$G^{k+n}$ has degree $(k+n)D$ and multiplicity $(k+n)M$ at each
point. By \cite[Theorem 1.1(a)]{HH02}, applied to the ideal $I$ of
the (reduced) scheme $\{p_1, \dots, p_r\}$, this implies the
existence of hypersurfaces of degree $t\le (k+n)D/M$ with
multiplicity at least $k+1$ at $p_1$, \dots, $p_r$. Now write
$\m'=(k+1, \dots, k+1)$; by \cite[Corollary 1.2]{AH00}, there is
$s_k(n)$ such that if $r\ge s_k(n)$ then $\alpha_{\m'}(p_1, \dots,
p_r) \ge (k+1)\sqrt[n]{r}$ (again, because the points are
general). Therefore we get $(k+n)D/M \ge (k+1)\sqrt[n]{r}$ and
$$\frac{d}{\sum m_i}=\frac{D}{rM}\ge \frac{1}{\sqrt[n]{r^{n-1}}}\cdot{\frac{k+1}{k+n}},$$
as wanted.
\end{proof}

\bibliographystyle{amsplain}
\bibliography{Biblio}
\end{document}